\documentclass[12pt,reqno]{amsart}
\usepackage{color,mathrsfs,epsfig}
\usepackage{psfrag,graphicx}
\usepackage{amssymb}

\newcommand{\N}{\mathbb{N}}

\newcommand{\R}{\mathbb{R}}

\newcommand{\Leb}{{\mathscr L}}

\newcommand{\ee}{\varepsilon}
\renewcommand{\div}{{\rm div}\,}
\newcommand{\divtx}{{\rm div}_{t,x}\,}

\newcommand{\loc}{{\rm loc}}

\renewcommand{\det}{{\rm det}}

\newtheorem{theorem}{Theorem}[section]
\newtheorem{lemma}[theorem]{Lemma}

\theoremstyle{definition}
\newtheorem{definition}[theorem]{Definition}

\theoremstyle{remark}

\newtheorem{conjecture}[theorem]{Conjecture}

\numberwithin{equation}{section}

\begin{document}
\bibliographystyle{plain}
\title[Transport equation and conservation laws]{An overview on some results concerning the transport equation and its applications to conservation laws}
\author{Gianluca Crippa}
\address{G.C.: Dipartimento di Matematica,
Universit\`a degli Studi di Parma,
viale G.P.~Usberti 53/A (Campus),
43100 Parma, Italy}
\email{gianluca.crippa@unipr.it}
\author{Laura V.~Spinolo}
\address{L.V.S.: Centro De Giorgi, Collegio Puteano, Scuola Normale Superiore,
Piazza dei Ca\-va\-lie\-ri 3, 56126 Pisa, Italy}
\email{laura.spinolo@sns.it}

\begin{abstract}
We provide an informal overview on the theory of transport equations with non smooth velocity fields, and on some applications of this theory to the well-posedness of hyperbolic systems of conservation laws.
\end{abstract}

\maketitle

\section{Introduction}
\label{intro}

Our attempt in this note is to give an overview on the theory of transport and continuity equations
$$
\partial_t u + b \cdot \nabla u = 0  \quad \text{ and } \quad
\partial_t u + \div ( bu) = 0 \,, \qquad
\text{ where $b : [0,+\infty[ \times \R^d \to \R^d$,}
$$ 
and on its applications to some classes of nonlinear PDEs. Since various extensive surveys and lecture notes on this topic are available (for instance the interested reader could consult \cite{Amb:cetraro,ambcri:ex,del:notes,thesis}, where also many references to the literature are provided), we prefer here to follow a more informal approach, by starting from the applications and then by step by step addressing the related issues concerning transport and continuity equations.

We now describe the topics that shall be discussed in this note: the related references will be provided in the corresponding sections.

We start in Section \ref{s:claws} by going over some of the main aspects of the theory of conservation laws. For this class of PDEs, which typically model conservation of physical quantities, existence fails in the class of classical (smooth) solutions, but bounded distributional solutions are not unique. Admissibility conditions have been designed in the attempt of restoring uniqueness of distributional solutions: an example is the entropy dissipation condition~\eqref{e:ntropia}, namely the requirement that some nonlinear function of the solution must be dissipated. The underling idea is that the entropy criterion should select the physically relevant solution. General existence and uniqueness  results concerning conservation laws are available in the scalar case, and for one-dimensional systems under smallness assumptions on the total variation of the initial data. However, proving existence and uniqueness results valid for general systems in several space dimensions is still a major open problem.

We then turn our attention to one particular multi-dimensional system, the so-called Keyfitz and Kranzer system, see \eqref{e:kk}. The peculiarity of \eqref{e:kk} lies in the dependence of the nonlinearity on the modulus of the solution only, thus inducing a sort of radial symmetry. Formally, this system splits into the scalar conservation law \eqref{e:rho} and the transport equation \eqref{e:theta}, and this paves the way for the application of transport equation techniques to the study of \eqref{e:kk}.

However, the only regularity (with respect to the space variable) we can expect on the velocity field appearing in \eqref{e:theta} is $BV$ regularity. Hence, we are forced to discuss the behaviour of this transport equation in the case of non smooth velocity fields. We go over this research area in Section \ref{s:transp}, addressing some of the main points of the well-posedness theory in the cases the velocity field enjoys either Sobolev or $BV$ regularity. 

We finally motivate in Section \ref{s:kk} how these techniques lead to rigorous well-posedness results for the Keyfitz and Kranzer system \eqref{e:kk}: actually, one has to introduce some additional arguments which take care of the fact that no information on the divergence, but only on the compression rate of the velocity fields are typically expected in this context. We conclude the note by collecting some results, analogue in spirit but valid in one space dimension, regarding again the Keyfitz and Kranzer system, and also the chromatography system \eqref{e:chrom}.

\section{Systems of conservation laws}
\label{s:claws}
In this section we informally go over some results concerning systems of conservation laws
in the form
\begin{equation}
\label{e:claw:gen}
  \partial_t U (t,x) + \sum_{i=1}^d \partial_{x_i} F_i \big[ U(t, x) \big]
     = 0 \,.
\end{equation}
The unknown $U$ is a function
$$
     U(t,x) : [0,+\infty[ \times \R^d \to \R^k
$$
and all the functions $F_i$, $i=1, \dots , d$ are defined on an open set of $\mathbb \R^k$ and take values in $\R^k$. Classical references concerning these topics are the books by Serre~\cite{Serre:book} and Dafermos~\cite{Daf:book}. Also, the book by Bressan~\cite{Bre:book} focuses on the study of the Cauchy problem for systems in one space dimension, $d=1$, while the recent book by Benzoni-Gavage and Serre~\cite{serrebenzoni} is concerned with systems in several space dimensions.

Conservation laws like~\eqref{e:claw:gen} have several and diverse physical applications. In particular, the Euler equation can be written in the form~\eqref{e:claw:gen}. The Euler equation is obtained from the Navier-Stokes equation by formally setting the viscosity coefficients equal to zero and describes the motion of the particles in an inviscid fluid.   

One of the main difficulties in the study of~\eqref{e:claw:gen} is the possible loss of regularity of the solution. More precisely, let us first focus on the Cauchy problem obtained by coupling~\eqref{e:claw:gen} with the initial condition
\begin{equation}
\label{e:caudat}
    U(0, x) =  \bar U (x).
\end{equation}
As a matter of fact, even if the datum~\eqref{e:caudat} is very regular, in general there is no classical solution of~\eqref{e:claw:gen} defined on the whole time interval $t \in  [0, \, +\infty[$.
Here a classical solution is a function $U \in C^1 \big( [0,+\infty[ \times \R^d; \R^k \big)$ satisfying~\eqref{e:claw:gen}-\eqref{e:caudat} pointwise. Remarkably, examples of solutions starting from a datum $\bar U \in C^{\infty}$ and developing discontinuities in finite time are available even in the scalar, one-dimensional case $k=d=1$.

In the attempt of obtaining global existence results, it is thus natural to introduce a notion of \emph{weak solution}, interpreting \eqref{e:claw:gen} in the sense of distributions. However, in general a weak solution of the Cauchy problem~\eqref{e:claw:gen}-\eqref{e:caudat} is not unique: as for the development of discontinuities, examples of non uniqueness for a distributional solution of~\eqref{e:claw:gen} are available even in the scalar, one-dimensional case $k=d=1$. In the attempt at selecting a unique solution, several \emph{admissibility conditions} can be introduced (entropy inequalities, Lax and Liu conditions on shocks), often motivated by physical considerations: we refer to the book by Dafermos~\cite[Chapters 4 and 8]{Daf:book} for a detailed exposition.

Since we will need it in the following, we recall the definition of entropy here. Let ${\eta: \R^k \to \R}$, $Q: \R^k \to \R^d$ be two smooth functions, $Q= (Q_1, \dots, Q_d)$. Then $(\eta, Q)$ is an entropy-entropy flux pair for~\eqref{e:claw:gen} if 
$$
    \nabla \eta \cdot D F_i = \nabla Q_i \quad \text{for every $i = 1 ,\dots , d$}.
$$  
Here $D F_i$ denotes the Jacobian matrix of $F_i$ and $\nabla \eta$ the gradient of $\eta$. A distributional solution of~\eqref{e:claw:gen} is entropy admissible if the inequality 
\begin{equation}
\label{e:ntropia}
    \partial_t \big[ \eta (U) \big]+ \div \big[ Q(U) \big] \leq 0
\end{equation}
holds in the sense of distribution for every $(\eta, Q)$ entropy-entropy flux pair such that $\eta$ is convex.  

In the scalar case $k=1$, the fundamental work by Kru{\v{z}}kov~\cite{kru} established existence and uniqueness of global in time, entropy admissible solutions of the Cauchy problem~\eqref{e:claw:gen}-\eqref{e:caudat} under the assumption that the initial datum $ \bar U$ is bounded, $\bar U \in L^{\infty} (\R)$. Also, Kru{\v{z}}kov~\cite{kru} showed that the $BV$ regularity is propagated, namely  $U(t, \cdot) \in BV (\R)$ for every $ t >0$ if $\bar U \in BV (\R)$. 

Let us now consider the one-dimensional case, $d=1$, $k \ge 1$: in this case,~\eqref{e:claw:gen} becames
\begin{equation}
\label{e:claw:oned}
           \partial_t U + \partial_x \big[ F(U) \big] =0, \quad \text{where $U \in \R^k$}. 
\end{equation}
The existence of a global in time, entropy admissible solution of the Cauchy problem~\eqref{e:claw:oned}-\eqref{e:caudat} was established by Glimm~\cite{glimm}. Uniqueness results were obtained in a series of papers by Bressan and several collaborators and are discussed in Bressan's book~\cite{Bre:book}. All these results concerning the Cauchy problem~\eqref{e:claw:oned}-\eqref{e:caudat} hold under the assumptions that the total variation of the initial datum $\bar U$ is small enough and that the system satisfies an hypothesis of so-called \emph{strict hyperbolicity} (see for example Dafermos~\cite[Sections 3.1 and 9.6]{Daf:book} for the exact definition and for a related discussion). If this is not the case, then either existence or uniqueness of global in time, entropy admissible solutions may fail (see e.g.~the examples discussed in Jenssen~\cite{Jenssen} and Panov~\cite{Panov:00}).
  
Finally, we consider the general case of a systems of conservation laws in several space dimensions~\eqref{e:claw:oned}, $d >1, k >1$. Presently, no general result is available concerning either existence or uniqueness of global in time distributional solutions. Tackling these issues in big generality is certainly regarded as one of the most challenging open problems in this research area.  

For specific classes of systems, however, well-posedess results are available. In the following, we focus on cases when these results are obtained by employing transport equation techniques. 

The Keyfitz and Kranzer system was introduced in~\cite{KK} and takes the form
\begin{equation}
\label{e:kk}
\left\{
\begin{array}{lll}
         \partial_t U +\displaystyle{ \sum_{\alpha = 1}^d \frac{\partial }{ \partial x_{\alpha}}
         \Big( f_{\alpha} ( | U |) U \Big) = 0 } \\ \\
         U(0, x) =  \bar U (x)\,.
\end{array}
\right.
\end{equation}
Here, $f_{\alpha} : \R \to \R$
is a smooth function for each $\alpha = 1, \dots , d$. By setting $\rho = |U|$ and 
$$
     \theta = (\theta_1, \dots , \theta_k) = U / |U|,
$$ system~\eqref{e:kk} \emph{formally} splits as the coupling between the scalar conservation law
\begin{equation}
\label{e:rho}
\left\{
\begin{array}{lll}
         \partial_t \rho +
         \div \big( f (\rho) \rho \big) =0  \\ \\
         \rho (0, x) = |  \bar U|(x)\,,
\end{array}
\right.
\end{equation}
and the transport equations 
\begin{equation}\label{e:theta}
\left\{
\begin{array}{lll}
    \partial_t \theta_i  +   f (\rho)  \cdot \nabla \theta_i  =0  \\ \\
    \theta_i (0,x) =  \bar U (x) / | \bar U (x) |
    \end{array}
\right.
\qquad \qquad i=1, \dots, k \,.
\end{equation}
Here, $f= (f_1, \dots, f_d)$. Some remarks are here in order: equation~\eqref{e:rho} is a scalar conservation law in several space dimensions and does not depend on $\theta$. We can thus apply Kru{\v{z}}kov's~\cite{kru} theory to obtain existence and uniqueness of global in time, entropy admissible, distributional solution. Then, we can plug the function $\rho$ in the second equation, obtaining a transport equation depending on $\theta$ only. However, even if the initial datum $\rho (0, x) = |\bar U (x)|$ is extremely regular, the best we can hope is that the $BV$ regularity is propagated, i.e.~that $\rho(t, \cdot) \in BV (\R)$ for every $ t >0$ if $\rho(0, \cdot) \in BV (\R)$. Indeed, as mentioned before, a solution having very regular initial data can in general develop discontinuities in finite time. 

In this way, we find that the analysis of~\eqref{e:rho}-\eqref{e:theta} has close links to the study of continuity and transport equations with weakly differentiable coefficients: an informal overview on these issues is provided in Section~\ref{s:transp}, while in Section~\ref{s:kk} we come back to the Keyfitz and Kranzer system~\eqref{e:kk} and we discuss some well-posedness results obtained by pursing the approach we sketched before.

\section{Transport equation} \label{s:transp}

In this section we are concerned with the well-posedness of the {\em transport equation}
\begin{equation}\label{e:cp2d}
\begin{cases} 
\partial_t u(t,x) + b(t,x) \cdot \nabla u(t,x) = 0 \\
u(0,x) = \bar{u}(x) 
\end{cases}
\end{equation}
and of the related {\em continuity equation}
\begin{equation}\label{e:cp2dcont}
\begin{cases} 
\partial_t u(t,x) + \div \big( b(t,x) u(t,x) \big) = 0 \\
u(0,x) = \bar{u}(x) \,,
\end{cases}
\end{equation}
in which $b : [0,+\infty[ \times \R^d \to \R^d$ is a vector field, $\bar u \in L^\infty(\R^d)$ is the initial datum, and the unknown $u$ belongs to $L^\infty([0,+\infty[ \times \R^d)$.

The well-poseness theory for these equations is classical and well understood in the case when $b$ is sufficiently smooth (at least Lipschitz with respect to the spatial variable, uniformly with respect to the time), and is strongly based on the so-called theory of characteristics, i.e.~on the connection between \eqref{e:cp2d} and \eqref{e:cp2dcont} and the {\em ordinary differential equation}
\begin{equation}\label{e:odebas}
\begin{cases}
\dot \gamma (t) = b(t,\gamma(t)) \\
\gamma(0) = x \,.
\end{cases}
\end{equation}
However, in many applications motived by physical models, non-smooth vector fields show up as velocity fields of transport or continuity equations. In addition to the case of the Keyfitz and Kranzer system~\eqref{e:kk} described at the end of Section~\ref{s:claws}, the theory of non-smooth transport equations has interesting applications to the study of the Boltzmann equation and of the Vlasov-Poisson equation: these applications were investigated by DiPerna and Lions in~\cite{boltzmann} and~\cite{vlasov} respectively.

This motivates the great interest arisen in the study of \eqref{e:cp2d} and \eqref{e:cp2dcont} when $b$ is only in some classes of weak differentiability. The first seminal result in this direction is due to DiPerna and Lions \cite{DiPLi}, who deal with the case of Sobolev regularity, under boundedness assumptions on the spatial divergence of the vector field. In \cite{DiPLi} the notion of renormalized solution to the transport (or continuity) equation is introduced. Let us explain heuristically the motivation for such a notion. Assume for the moment that the vector field $b$ is compactly supported and that the spatial divergence
\begin{equation}\label{e:assdiv}
\div b \; \in \; L^1_t (L^\infty_x) = L^1 \big( [0,+\infty[_t ; L^\infty(\R^d_x) \big) \,.
\end{equation}
We multiply our equation
\begin{equation}\label{e:pre}
\partial_t u + b \cdot \nabla u = 0
\end{equation}
by $2u$, and by formally applying the chain-rule we deduce
\begin{equation}\label{e:post}
\partial_t \big( u^2 \big) + b \cdot \nabla \big( u^2 \big) =0 \,.
\end{equation}
We now integrate over the space $\R^d$ for every fixed valued $t$ of the time, obtaining
\begin{equation}\label{e:needdiv}
\begin{split}
\frac{d}{dt} \int_{\R^d} u^2(t,x) \, dx &= \int_{\R^d} \div b(t,x) u^2(t,x) \, dx \\
&\leq \| \div b (t,\cdot) \|_{L^\infty(\R^d)} \int_{\R^d} u^2(t,x) \, dx \,,
\end{split}
\end{equation}
thanks to the divergence theorem. Now, a simple application of Gronwall's lemma implies that, if the $L^2(\R^d)$ norm of the solution vanishes at the initial time, then it vanishes for all time. By the linearity of \eqref{e:cp2d} we obtain uniqueness for the Cauchy problem.

In the above argument, the delicate point is hidden in the passage from \eqref{e:pre} to \eqref{e:post}. Indeed, since no regularity beyond boundedness is assumed on the solution $u$, the application of the chain-rule formula is not justified. This led DiPerna and Lions \cite{DiPLi} to define {\em renormalized solution} of \eqref{e:cp2d} as bounded distributional solutions $u$ for which 
\begin{equation}\label{e:renfor}
\partial_t \big( \beta(u) \big) + b \cdot \nabla \big( \beta(u) \big) = 0
\end{equation}
holds for any function $\beta \in C^1(\R;\R)$. In some sense, the validity of the ``chain-rule formula'' is assumed ``by definition'' for renormalized solutions.

By relying on an argument similar to the one we gave before, one can show that uniqueness results are satisfied by the solutions of the transport equation~\eqref{e:cp2d} and the continuity equation~\eqref{e:cp2dcont}, provided that the vector field $b$ satisfies the following property: all bounded distributional solutions of~\eqref{e:pre} are renormalized and satisfy some strong continuity requirement.

The core of the argument in DiPerna and Lions~\cite{DiPLi} is thus the proof of the validity of \eqref{e:renfor} for Sobolev vector fields, and this is achieved via a regularization procedure combined with a control on the convergence of the error term that appears. This gives well-posedness for \eqref{e:cp2d} (and similarly for \eqref{e:cp2dcont}) assuming that $b$ is a bounded vector field such that
\begin{equation}
\label{e:dipl}
b \in L^1_t(W^{1,1}_x) \qquad \text{ and } \qquad \div b \in L^1_t(L^\infty_x) \,.
\end{equation}
More recently, this result has been extended by Ambrosio \cite{Amb:trabv} to the case of bounded vector fields with bounded variation with respect to the space variable (i.e.~such that the space distributional derivative is a locally finite measure), that is one replaces~\eqref{e:dipl} with
$$
b \in L^1_t(BV_x) \qquad \text{ and } \qquad \div b \in L^1_t(L^\infty_x) \,.
$$
In both results, the need of controlling the spatial divergence of $b$ essentially comes from computations analogue to \eqref{e:needdiv}. The hard point of the proof in \cite{Amb:trabv} is to control the convergence of the error term in the regularization procedure in this weaker context. The argument is based on more refined arguments of geometric measure theory, in particular on some fine properties of $BV$ functions. For a general survey on this topic, we suggest for instance \cite{Amb:cetraro,ambcri:ex,thesis}.

The possibility of extending the uniqueness result by removing the assumption 
$b \in  L^1_t( BV_x)$ is ruled out by an intriguing counterexample due to Depauw~\cite{depauw:ex}: he exhibits a bounded, divergence free vector field such that the Cauchy problem
$$
    \left\{
    \begin{array}{ll}
   \partial_t u + \div (b u ) =0 \\
   u(0, x) \equiv 0
   \end{array}
   \right.
$$
admits a non trivial solution. The vector field in~\cite{depauw:ex} satisfies  $b \in L^1_{\loc} \big( ]0; 1]; BV( \R^2; \R^2) \big)$, but the $BV$ norm is not integrable up to time $t=0$. 

As a side remark, we mention that a result in a recent work~\cite{ACFS} by Ambrosio, Figalli and the authors establishes well-posedness for continuity equations where the vector field $b \in BV_{\loc}  ( ] 0, + \infty[ \times \R^d; \R^d)$ 
and has a possible blow up of the $BV$ norm at $t=0$. This does not contradict Depauw's counterexample because this result is obtained by imposing on the solution an additional regularity assumption, namely that the map $t \mapsto u(t, \cdot)$ is continuous in $L^1_{\loc} (\R^d)$ endowed with the strong topology. 

\section{Applications to systems of conservation laws}
\label{s:kk}
\subsection{Nearly incompressible vector fields}
Let us go back to the example of the Keyfitz and Kranzer system~\eqref{e:kk}: as pointed out at the end of Section~\ref{s:claws}, Kru{\v{z}}kov's~\cite{kru} theory of entropy solutions ensures  that the vector field 
$$ 
    b(t, x) = f (\rho) (t, x)
$$
appearing in the transport equation~\eqref{e:theta} has $BV$ regularity with respect to the space variable. Thus, the theory described in Section \ref{s:transp} is concerned with the right regularity framework, i.e.~with vector fields having bounded variation. However, the additional assumption of boundedness of the divergence of the vector field may be not satisfied, and this means that the theory does not apply straightforwardly to the analysis of the Keyfitz and Kranzer system~\eqref{e:kk}. Remember that boundedness of the divergence of the vector field is apparently an unavoidable assumption if one tries to carry on computations like in \eqref{e:needdiv}.

We now illustrate some assumptions one can impose on the vector field $b$ in~\eqref{e:cp2d} that are more natural in view of the applications to the Keyfitz and Kranzer system~\eqref{e:kk}: in this case, we have an additional information, coming from \eqref{e:rho}, namely that the function $\rho$ is a $BV$ solution of the continuity equation with vector field $f(\rho)$. Also, Kru{\v{z}}kov's theory~\cite{kru} guarantees that the scalar conservation law~\eqref{e:rho} satisfies a comparison principle, namely 
$$
    0 < \frac{1}{C} \leq \rho (t, x) \leq C \quad \text{for every $t \ge 0$ and for a.e.~$x \in \R$} 
$$
if $
    1 / C \leq \rho (0, x) \leq C  
$
for a.e.~$x \in \R$. Note that requiring that $\rho (0, x)$ is non negative is natural here, since $\rho (0, x)$ is the modulus of the initial datum. 

Let us consider the $(d+1)$-dimensional vector field in $\R^+_t \times \R^d_x$ defined by
$$
B(t,x) = \big( \rho(t,x) , f(\rho(t,x)) \rho(t,x) \big) \in \R \times \R^d\,.
$$
Thus, \eqref{e:rho} simply means that $B(t,x)$ is divergence free in space-time, i.e.
$$
\divtx B(t,x) = 0 \,.
$$

We add a fictitious ``time variable'' $s \in [0,+\infty[$, and by setting 
$$
    \tilde \theta_i (s,t,x) = \theta_i (t, x)
$$  we realize that the equation written at the first line of~\eqref{e:kk} 
can be rewritten as
\begin{equation}\label{e:lifttrans}
\partial_s \tilde \theta_i (s,t,x) + \divtx \big( B(t,x) \tilde \theta_i (s,t,x) \big) = 0 \, \quad \text{for $i=1 ,\dots, k$}.
\end{equation}
Here we are exploiting the fact that
$$
    \theta (t, x) = \big( \theta_i ,\dots, \theta_k \big) (t, x) = U / |U|. 
$$
Thus, we can now apply the results of Section \ref{s:transp} to the ``lifted'' equation \eqref{e:lifttrans}, since we are now dealing with a $BV$ divergence free vector field (precisely thanks to the $BV$ regularity of $\rho$). This means, in particular, that $\tilde \theta_i (s,t,x)$ is a renormalized solution of \eqref{e:lifttrans}, that is
$$
\partial_s \big[ \beta \big(\tilde \theta_i (s,t,x) \big)\big] + 
\divtx \big[ B(t,x) \beta \big( \tilde \theta_i (s,t,x) \big) \big] = 0
$$
for all $\beta \in C^1(\R;\R)$. But, $\tilde \theta_i$ not depending on $s$, this can be rewritten as
\begin{equation}\label{e:renincompr}
\partial_t \big[ \rho(t,x) \beta\big(\theta_i(t,x)\big) \big] + 
\div \big[ b(t,x) \rho(t,x) \beta\big(\theta_i(t,x)\big) \big] = 0\,,
\end{equation}
which is reminiscent of the renormalization condition for the vector field $b(t,x)$, apart from the term $\rho (t,x)$ now appearing. 

The previous discussion motivates the following definition:
\begin{definition}
\label{d:nivf}
We say that a bounded vector field $b$ is {\em nearly incompressible} if there exist a function $\rho$ and a constant $C>0$ such that 
\begin{equation}\label{e:costC}
0 < \frac{1}{C} \leq \rho(t,x) \leq C < +\infty \qquad \text{for $\Leb^{d+1}$-a.e.~$(t,x) \in \, ]0, + \infty[ \times \R^d$}
\end{equation} 
and
\begin{equation}\label{e:near} 
\partial_t \rho + \div ( b \rho ) = 0 \,.
\end{equation}
\label{d:incompr}\end{definition}

The renormalization property, in this context, is the requirement of the validity of \eqref{e:renincompr} for every bounded function $\theta_i$ solving
$$
    \partial_t \big[ \rho \theta_i \big] + \div \big[ b \rho \theta_i \big] = 0
$$

One can show that uniqueness results hold for nearly incompressible vector fields with the renormalization property: we refer to the notes by De Lellis~\cite{del:notes} for a complete and systematic treatment of the renormalization theory for nearly incompressible vector fields and for various related references. The following lemma is one of the main results of this theory, and can be proved by relying on the ``fictitious variable" argument described above. 
\begin{lemma}[{\cite[Lemma 5.10]{del:notes}}]\label{l:lemmanote} Let 
$$
b \in L^\infty  \cap BV_\loc \big( [0,+\infty[ \times \R^d ; \R^d \big)
$$
and 
$$
\rho \in L^\infty  \cap BV_\loc \big( [0,+\infty[ \times \R^d \big)\,,
\qquad \text{ with $\rho(0,\cdot) \in BV_\loc(\R^d)$}
$$
be such that
$$
\partial_t \rho + \div \big( b \rho \big) = 0 \,.
$$
Let $u \in L^\infty \big( [0,+\infty[ \times \R^d \big)$ be a solution of
$$
\partial_t \big( \rho u \big) + \div \big( b \rho u \big) = 0 \,.
$$ 
Then $u$ is renormalized, that is for every $\beta \in C^1(\R;\R)$ there holds
$$
\partial_t \big( \rho \beta(u) \big) + \div \big( b \rho \beta(u) \big) = 0 \,.
$$ 
\end{lemma} 
 
A major conjecture concerns the possibility of removing the assumption of $BV$ regularity on $\rho$ in Lemma \ref{l:lemmanote}:
\begin{conjecture}\label{c:renconj}
Any bounded nearly incompressible vector field $ b \in BV ( [0, + \infty[ \times \R^d; \R^d)$ has the renormalization property. 
\end{conjecture}

This conjecture is also important in view of a compactness conjecture advanced by Bressan in \cite{bre:conj}, again in connection to the question of existence of solutions to the Keyfitz and Kranzer system. Indeed, Ambrosio, Bouchut and De Lellis~\cite{ABDL} proved that a positive answer to Conjecture \ref{c:renconj} would imply a positive answer to the compactness conjecture, which reads as follows:

\begin{conjecture}[Bressan's compactness conjecture]\label{c:bressan2}
Let $b_k: [0,+\infty[ \times \R^d\to \R^d$, $k\in \N$, be a sequence of smooth vector fields
and denote by $X_k$ the (classical) flows associated to them.
Assume that $\|b_k\|_\infty + \|\nabla b_k\|_{L^1}$ 
is uniformly bounded and that
the flows $X_k$ are uniformly nearly incompressible, i.e.~that
\begin{equation}\label{e:incompres102}
\frac{1}{C} \;\leq\; \det ( \nabla_x X_k (t,x)) \;\leq\; C\, 
\qquad \text{for some constant $C>0$}.
\end{equation}
Then the sequence
$\{X_k\}$ is strongly precompact in $L^1_\loc([0,+\infty[ \times \R^d)$.
\end{conjecture}

Both conjectures are presently open, apart from some particular cases of regularity. Some recent advances are presented in  Ambrosio, De Lellis and Mal\'y~\cite{Amb-De-Ma}.

\subsection{Keyfitz and Kranzer system in several space dimensions}
By relying on Lem\-ma~\ref{l:lemmanote} and on the uniqueness results for nearly incompressible vector fields with the normalization property, Ambrosio and De Lellis~\cite{ambdel} and Ambrosio, Bouchut and De Lellis~~\cite{ABDL} proved well-posedness results for the Keyfitz and Kranzer system~\eqref{e:kk}. They followed the approach sketched at the end of Section~\ref{s:claws}, which was inspired by considerations in Bressan~\cite{Bre:illposed}. 

Before stating the well-posedness result in~\cite{ABDL}, we have to introduce a preliminary definition. Let $U$ be a locally bounded function solving~\eqref{e:kk} in the sense of distributions. Then $U$ is a \emph{renormalized entropy admissible solution} if $\rho = |U|$ is an entropy admissible solution of~\eqref{e:rho} satisfying ${\lim_{t \to 0^+} \rho (t, \cdot) =|  \bar U|}$ in the strong topology of $L^1_\loc(\R^d)$.
By relying on a classification of the entropies for~\eqref{e:kk} due to Panov~\cite{Panov:00} and Frid~\cite{Frid}, under quite general assumptions on the function $f$ any renormalized entropy solution is indeed an entropy admissible solution of~\eqref{e:kk} (for a proof of this implication see for example the notes by De Lellis~\cite[Proposition 5.7]{del:notes}). Conversely, in general entropy admissible solutions of~\eqref{e:kk} are not renormalized and they are not unique (for an example see e.g.~Bressan~\cite{Bre:illposed}).

Existence, uniqueness and stability of renormalized entropy admissible solutions are established in the following theorem.
\begin{theorem}[{\cite[Theorem 2.6]{ABDL}}]
         Assume that the function $f = (f_1 ,\dots, f_d)$ is locally Lipschitz continuous and that the modulus of 
         the initial datum   satisfies $|U| \in L^{\infty} (\R^d) \cap BV_{\loc} (\R^d)$. Then there exists a unique renormalized entropy admissible solution of~\eqref{e:kk}. Also, assume that $\bar U_n$ is a sequence of initial data satisfying the following assumptions:
         \begin{enumerate}
         \item $ \| \bar U_n  \|_{L^{\infty}} \leq C$ for some constant $C>0$;
         \item for every bounded open set $\Omega \subset \R^d$, $\bar \rho_n = |U_n|$ 
         satisfies $\|\bar \rho_n \|_{BV (\Omega)} \leq C(\Omega)$;
                  \item $\bar U_n \to \bar U$ strongly in $L^1_{\loc} (\R^d, \R^k)$. 
         \end{enumerate}
         Then $U_n \to U$ strongly in $L^1_{\loc} ( [0, + \infty[ \times \R^d, \R^k).$  
\end{theorem}

As a final remark, we point out that a counterexample due to Crippa and De Lellis~\cite{CDL} indicates that strategies exploiting transport equation techniques aimed to show well-posedness for hyperbolic systems are not expected to work for {\em general} systems of conservation laws. More precisely, the analysis in~\cite{CDL} shows that (in dimension greater than or equal to $3$) there exists no functional space containing $BV$ which is closed by iteration of transport equations and which embeds compactly in $L^1_\loc$. This suggests that approximation schemes valid for general conservation laws cannot be constructed by relying on transport equation techniques.  
\subsection{Transport equations techniques applied to some systems in one space dimension}
If one restricts to one space dimension $d=1$, well-posedness results for the continuity equation
$$
    \partial_t u + \partial_x (bu) =0
$$
are available under much weaker assumptions than those discussed in the previous sections. In particular, Panov~\cite{Panov:99, panov:gensol} proved existence and uniqueness results for bounded vector fields $b \in L^{\infty} ([0, + \infty [ \times \R)$ satisfying a near incompressibility condition similar to the one 
introduced in Definition 4.1. See also Serre~\cite{serre:onde}. These results are then applied to the analysis of the one-dimensional Keyfitz and Kranzer system
$$
   \partial_t U + \partial_x \big[ f (|U|) U \big] =0
$$
obtaining existence and uniqueness of renormalized, entropy admissible solutions with merely bounded initial data (the presence of linear source terms is also contemplated in~\cite{panov:gensol}). See also Panov~\cite{Panov:00}, where the vanishing viscosity approximation 
$$
   \partial_t U^{\ee} + \partial_x \big[ f (|U^{\ee}|) U^{\ee} \big] = \ee U^{\ee}_{xx},
$$
is investigated, and the convergence of the family $U^{\ee}$ as $\ee \to 0^+$ is established. 

As an other example of applications of transport equation techniques to a system of conservation laws in one dimension we mention the so-called chromatography system 
\begin{equation}\label{e:chrom}
\begin{cases}
\partial_t u_1 + \partial_x \left( \displaystyle \frac{u_1}{1+u_1+\cdots + u_k} \right) = 0 \\ \\
\qquad\qquad\vdots \\ \\
\partial_t u_k + \partial_x \left( \displaystyle \frac{u_k}{1+u_1+ \cdots + u_k} \right) = 0\,.
\end{cases} 
\end{equation}
Here the unknowns $u_1, \dots , u_k$ are real-valued and one is usually interesting in finding a solution satisfying 
$u_1 \ge 0, \dots , u_k \ge 0$. By setting 
$$
    v = u_1 + \cdots + u_k \quad  \text{and} \quad  w_i = u_i \quad \text{for $i=2, \dots, k$}
$$
one obtains that the system writes as the coupling between the scalar conservation law 
$$
\partial_t v + \partial_x \left( \displaystyle \frac{v}{1+v} \right) = 0 
$$
and the continuity equations
\begin{equation}
\label{e:w}
\partial_t w_i + \partial_x \left( \displaystyle \frac{w_i}{1+v} \right) = 0 \quad \text{for $i=2, \dots, k$}.
\end{equation}
As in~\eqref{e:theta}, the vector field $b (t, x) = 1 / (1 + v)$ in~\eqref{e:w} depends on $v$ and hence one expects low regularity. By following this strategy, well-posedness results for the chromatography system are obtained in~\cite{ACFS}. See also the previous works by Panov~\cite[Remark 4 page 140]{Panov:00} and Bressan and Shen~\cite{bressanshen} for similar approaches. 
\bibliography{bibliochroma}

\end{document}